\documentclass[A4, 11pt]{article}

\usepackage{geometry}

\usepackage{tikz}
\tikzset{every picture/.style={scale=0.5}}

\usepackage{graphicx}
\usepackage{float}
\usepackage{color}
\usepackage{enumitem}
\usepackage{subcaption}
\usepackage{caption}

\usepackage{authblk}

\usepackage{amssymb,amsmath,amsthm}

\usepackage{mathtools}
\usepackage[linktoc=all,hidelinks]{hyperref}

\DeclareMathOperator{\ord}{ord}
\DeclareMathOperator{\supp}{supp}

\newcommand{\bra}[1]{\langle#1\rangle}

\newcommand{\G}[1]{\mathbb Z_{#1}\times\mathbb Z_{#1}}   
\newcommand{\F}[1]{\widehat{\mathbf 1}_{#1}}
\newcommand{\Z}[1]{Z(\F{#1})}

\newtheorem{Thm}{Theorem}
\newtheorem{Lm}{Lemma}
\newtheorem{Cor}{Corollary}

\title{An Equi-distribution Equality of Tiles}
\author{Weiqi Zhou\thanks{zwq@xzit.edu.cn}}
\affil{\small School of Mathematics and Statistics, Xuzhou University of Technology \\  {\footnotesize Lishui Road 2, Yunlong District, Xuzhou, Jiangsu Province, China 221111}}
\date{}							

\begin{document}
\maketitle
\begin{abstract}
We prove that if the zero set of the Fourier transform of $A\subseteq\G{n}$ contains an element of prime power order, then there is an equi-distribution relation in subsets of $A$ with respect to certain hyperplanes. With this we further show that if $A$ is a tiling complement of the subgroup generated by $(p,0)$ and $(0,p)$ in $\G{p^m}$, then the zero set of its Fourier transform is disjoint with the orthogonal rotation of $A$. These results are motivated by a casual observation in $\G{p^2}$. \\

{\noindent
}
\end{abstract}
 
\section{Motivation}
Let $\mathbb Z_n$ be the additive cyclic group with order $n$. A subset $A$ in $\G{n}$ is called a \emph{tiling set} (or a \emph{tile}) if there exists some $B\subseteq\G{n}$ such that $\{A+b\}_{b\in B}$ forms a partition of $\G{n}$, $(A,B)$ is then said to be a \emph{tiling pair} in $\G{n}$ and we will write $\G{n}=A\oplus B$. Alternatively $A$ is called a \emph{Euclidean spectral set} if there exists some $S\subseteq\G{n}$ such that $\{e^{2\pi i\bra{a,s}/n}\}_{s\in S}$ (Here $\bra{a,s}$ is the usual Euclidean inner product) forms an orthogonal basis on $L^2(A)$ with respect to the counting measure. $(A,S)$ is then said to be a \emph{Euclidean spectral pair} in $\G{n}$. If $(A,B)$ is a tiling pair, then we will call them \emph{tiling complements} of each other in $\G{n}$. Alternatively if $(A,S)$ is a Euclidean spectral pair, then we will call $S$ a \emph{Euclidean spectrum} of $A$.

A fundamental question of the sampling theory, known as the Fuglede conjecture \cite{fuglede1974}, asks whether these two types of sets coincide. The conjecture has been disproved in $\mathbb R^d$ for $d\ge 3$ \cite{farkas2006, matolcsi2006, kolounzakis2006, matolcsi2005, tao2004}, while for $d=1,2$ it remains open. Their coincidence in finite Abelian groups are equally important, and are partial equivalent to the same conjecture in $\mathbb R^d$  \cite{dutkay2014, iosevich2013, wang1996}. On $p$-adic fields it is known to be true \cite{fan2016, fan2019}. 

So far there are quite a number of positive results in cyclic groups, in which the methods evolves. To keep the overall size of the paper short they are not to be fully listed here. For groups generated by two elements, the number of results (all positive) are limited and methods therein are mostly combinatorial \cite{fallon2022, Iosevich2017, shi2020, zhang2021}.

In the study of tiles in $\G{p^2}$,  the author observed an interesting phenomenon that if $(A,B)$ is a tiling pair with $|A|=|B|$, then sometimes $(A,B)$ is also a spectral pair. So why should we even consider taking $B$ as a spectrum of $A$ given that it is already a tiling complement of $A$? The initial motivation concerns the universal spectrum property: If $S$ is a spectrum for any tiling complement of $A$, then it is called a \emph{universal spectrum} (for tiling complements) of $A$. It is known that all tiling sets are spectral if and only if every tiling set admits a universal spectrum for its tiling complements \cite{farkas2006, kolounzakis2006, wang1997}. 

The motivation is to think reversely: Suppose that $A$ is a tiling complement of an order $n$ subgroup $H$ in $\G{n}$, and $B$ is a tiling complement of $A$ (it is possible for $B\neq H$ to hold, for example if $A=\{(0,0),(3,0),(0,3),(3,3)\}$, then $A$ is a tiling complement of the non-cyclic subgroup $\{(0,0),(2,0),(0,2),(2,2)\}$ in $\G{4}$, but it is, at thes same time, also a tiling complement of $\{(0,0), (1,2), (2,0), (3,2)\}$), and assuming the existence of a universal spectrum, then $B$ and $H$ should share this universal spectrum of $A$, a very natural choice of this spectrum would be $A$ itself since in many cases (though not always) $A$ and $H$ are actually spectra of each other (for instance, if $H=\mathbb Z_n\times\{0\}$, both $a_0,\ldots,a_{n-1}$ and $b_0,\ldots,b_{n-1}$ are permuted sequences of $0,\ldots,n-1$, then $A=\{(a_0,b_0),\ldots,(a_{n-1},b_{n-1})\}$ is both a tiling complement and a spectrum of $H$ at the same time, while $\{0\}\times\mathbb Z_n$ is only a tiling complement but not a spectrum of $H$). A significant case is in $\G{p^2}$ when $A=\{0,\ldots,p-1\}\times\{0,\ldots,p-1\}$ and $B$ is the (only) non-cyclic subgroup that consists of all order $p$ elements, then not only that they are tiling pair and spectral pair at the same time, but also the zeros sets of their respective Fourier transforms are disjoint with their own difference sets respectively, consequently both of them are actually universal spectra for their respective tiling complements. However, even though the assumption of  $A$ being a tiling complement of some non-cyclic subgroup is very structural, at the moment there is still very little we know about the tiling pair $(A,B)$. 

The purpose of this short note is to prove in Theorem \ref{ThmCount} that if the Fourier transform of $A\subseteq\G{n}$ has an element of prime power order in its zero set, then there is an equi-distribution relation in subsets of $A$ with respect to certain hyperplanes. With this we show further in Corollary \ref{CorDiff} that if $A$ is a tiling complement of the subgroup generated by $(p,0)$ and $(0,p)$ in $\G{p^m}$, then the zero set of its Fourier transform of $A$ is disjoint with the orthogonal rotation of $A$. Theorem \ref{ThmCount} is previously given by the author using the symplectic Fourier transform, in this paper all the derivation are carried out with the usual Fourier transform.

\section{Preliminaries} 
The \emph{difference set} of $A\subseteq\G{n}$ is
$$\Delta A=\{a-a': a,a'\in A, a\neq a'\}.$$
Suppose $A,B\subseteq\G{n}$ , and $a,a'\in A, b,b'\in B$, then
$$a+b=a'+b' \quad\Leftrightarrow\quad a-a'=b-b',$$
which means
$$A\oplus B \text{ is well defined} \quad\Leftrightarrow\quad \Delta A\cap \Delta B=\emptyset,$$
thus $(A,B)$ being a tiling pair in $\G{n}$ is equivalent to
$$|A|\cdot |B|=|\G{n}|, \quad\text{and}\quad \Delta A\cap \Delta B=\emptyset.$$

Let $f$ be a function on $\G{n}$, denote by  
$$\supp(f)=\{x\in\G{n}: f(x)\neq 0\}, \quad Z(f)=\{x\in\G{n}: f(x)=0\},$$
its support set and its zero set respectively. Its \emph{Fourier transform} is defined as
$$\widehat{f}(\xi)=\sum_{x\in \G{n}}f(x)e^{2\pi i\bra{x,\xi}/n}, \quad \xi\in\G{n}.$$
In particular, if $f=\mathbf 1_A$ is the characteristic function for some $A\subset\G{n}$, then we get
$$\F{A}(\xi)=\sum_{a\in A}e^{2\pi i\bra{a,\xi}/n}.$$

If $A$ is a finite multiset (i.e., a set that allows repetitive elements) on $\G{n}$ (i.e., each element of $A$ is a member of $\G{n}$), and the multiplicity of an element $a\in A$ (i.e., the number of copies of $a$ in $A$) is denoted by $m_a$, then the characteristic function on such a multiset $A$ can be defined and written as
$$\mathbf 1_A(x)=\begin{cases}m_x & x\in A, \\ 0 & x\notin A. \end{cases}$$
Its Fourier transform on $\mathbb Z_n^d$ is 
$$\F{A}(\xi)=\sum_{a\in A}m_a\cdot e^{2\pi i\bra{a,\xi}/n}.$$

If $A,B$ are finite multisets, and $A+B$ is the multiset obtained by enumerating $a\in A, b\in B$ (so $a$ appears $m_a$ times and $b$ appears $m_b$ times) and take the collection of all sums $a+b$, then straightforward computation shows
$$\F{A+B}=\F{A}\cdot \F{B}.$$

With these settings we easily get that when $A,B,S$ are usual sets (every element appear only once), then
\begin{align*}
(A,S) \text{ is a spectral pair on } \G{n} &\quad\Leftrightarrow\quad  \begin{cases}|A|=|S|, \\ \Delta S\subseteq \Z{A}, \end{cases} \\
(A,B) \text{ is a tiling pair on } \G{n} &\quad\Leftrightarrow\quad \begin{cases}|A|\cdot|B|=n^2, \\ \Delta\G{n}\subseteq Z(\F{A}\cdot\F{B}),\end{cases} 
\end{align*}
 
\section{Main Results}
From now on we use $\bra{x}$ for the cyclic subgroup generated by $x$, and $\ord(x)$ for the order of $x$ in $\G{n}$.

Given $A\subseteq \G{n}$, define its orthogonal set to be 
$$A^{\perp}=\{x\in\G{n}: \bra{a,x}=0\}.$$
It is easy to verify that $A^{\perp}$ is a subgroup of $\G{n}$, and 
$$x\in A^{\perp} \Rightarrow (\bra{x}\setminus\{0\}) \subseteq A^{\perp}.$$

Given $A\subseteq \G{n}$ and $k\in\mathbb Z_n$ we set
$$V_{x}^{(k)}=\{x\in\G{n}: \bra{a,x}=k \pmod n \}.$$
In particular, it follows immediately from the definitions that 
$$V_{x}^{(0)}=\bra{x}^{\perp}.$$

\begin{Thm} \label{ThmCount}
Let $p$ be a prime number, $x\in\G{n}$ with $\ord(x)=p^m$ ($m\in\mathbb N$) and $A\subseteq\G{n}$, if $x\in\Z{A}$, then we have
\begin{equation} \label{EqCount}
|A\cap\bra{px}^{\perp}|=|A\cap\bra{x}^{\perp}|\cdot p,
\end{equation}
and
\begin{equation} \label{EqEqui}
|A\cap V_{x}^{(0)}|=|A\cap V_{x}^{(n/p)}|=\ldots=|A\cap V_{x}^{(n(p-1)/p)}|.
\end{equation}
\end{Thm}

\begin{proof}
Since $\ord(x)=p^m$, there exists a unique $h\in\G{n}$ that satisfies $\ord(h)=n$ and 
$$x=\frac{n}{p^m}\cdot h.$$

Setting $\omega=e^{2\pi i/p^m}$ we get
\begin{equation}  \label{EqCountd3}
\F{A}(x)=\sum_{a\in A}e^{2\pi i\bra{a,x}/n}=\sum_{a\in A}e^{2\pi i\bra{a,h}/p^m}=\sum_{a\in A}\omega^{\bra{a,h}}.
\end{equation}

We shall view the right-hand side of \eqref{EqCountd3} as the polynomial $\sum_{a\in A}z^{\bra{a,h}}\in\mathbb Z[z]$ evaluated at $\omega$, then $\F{A}(x)=0$ implies that the polynomial is divisible by the $p^m$-th cyclotomic polynomial $1+z^{p^{m-1}}+\ldots+z^{(p-1)p^{m-1}}$, therefore we get 
$$\sum_{a\in A}\omega^{\bra{a,h}}=\left(\sum_{j=0}^{p-1}\omega^{j\cdot p^{m-1}}\right) \cdot \left(\sum_{k=0}^{p^{m-1}-1}c_k\omega^k\right),$$
with $c_k\in\mathbb Z$. Expand it further we obtain

\begin{equation} \label{EqExpand}
\left(\sum_{j=0}^{p-1}\omega^{j\cdot p^{m-1}}\right) \cdot \left(\sum_{k=0}^{p^{m-1}-1}c_k\omega^k\right)=c_0(1+\omega^{p^{m-1}}+\ldots+\omega^{(p-1)p^{m-1}})+g(\omega),
\end{equation}
where
$$g(\omega)=\left(\sum_{j=0}^{p-1}\omega^{jp^{m-1}}\right) \cdot \left(\sum_{k=1}^{p^{m-1}-1}c_k\omega_k\right),$$
collects all terms of form $c_k\omega^k$ with $k$ being a number that is not divisible by $p^{m-1}$. 

Compare the part $c_0(1+\omega^{p^{m-1}}+\ldots+\omega^{(p-1)p^{m-1}})$ in \eqref{EqExpand} with \eqref{EqCountd3} we see that for each $k\in\mathbb Z_p$, the term $ c_0\omega^{kp^{m-1}}$ must be obtained from those elements $a$ of $A$ that satisfies
$$\bra{a,h}=kp^{m-1},$$
which is equivalent to
$$\bra{a,x}=\frac{n}{p^m}\cdot\bra{a,h}=\frac{nk}{p},$$
i.e., $a\in V_x^{(nk/p)}$. 
Counting the number of terms we get 
$$c_0=|A\cap\bra{x}^{\perp_s}|=|A\cap V_{x}^{(0)}|=|A\cap V_{x}^{(n/p)}|=\ldots=|A\cap V_{x}^{(n(p-1)/p)}|,$$
which is \eqref{EqEqui}. 

For any $a\in V_x^{(nk/p)}$, it is also clear that 
$$\bra{a,px}=p\bra{a,x}=0\pmod n,$$
therefore
$$|A\cap\bra{px}^{\perp}|=\sum_{k\in\mathbb Z_p}|A\cap V_{x}^{(nk/p)}|=c_0p=|A\cap\bra{x}^{\perp_s}|\cdot p,$$
which establishes \eqref{EqCount}. 
\end{proof}

In particular, \eqref{EqEqui} described an equi-distribution phenomenon as described in \cite{shi2020} and references therein.

Given $x=(x_1,x_2)\in\G{n}$, define the orthogonal rotation of $x$ to be
$$x^{\Lsh}=(-x_2,x_1),$$
then it is obvious that
$$\ord(x)=\ord(x^{\Lsh}), \quad x^{\Lsh}\in\bra{x}^{\perp}.$$
Clearly $x\mapsto x^{\Lsh}$ is a bijection on $\G{n}$, its inverse is 
$$x^{\Rsh}=(x_2,-x_1).$$

For $A\subseteq\G{n}$, define its orthogonally rotated set to be the set obtained by rotating every element of $A$ orthogonally, i.e.,
$$A^{\Lsh}=\{a^{\Lsh}: a\in A\}.$$

In contrast to the intuition one may have acquired from the Euclidean spaces, here it is possible for $A=A^{\Lsh}=A^{\Rsh}$ to hold (i.e., a set can be orthogonal to itself). And unlike the usual ``orthogonal decomposition'' in Euclidean spaces, here it is also possible for a maximal cyclic subgroup to have no orthogonal complement. One example would be the subgroup generate be $(1,18)$ in $\G{25}$, it is easy to see that $(-18,1)=7\cdot (1,18) \pmod{25}$, so the subgroup is orthogonal to itself and one can also verify (for example, by using a computer) that it has no orthogonal complement. In general, this can happen in $\G{p^m}$ when $p \bmod 4=1$, since then there would be some $k$ satisfying 
$$k^2=-1 \pmod{p^m},$$
and $(-k,1)=k^{-1}\cdot (1,k)$ where $k^{-1}$ is the multiplicative inverse of $k$ modulo $p^m$.

\begin{Lm} \label{LmSelf}
Let $p$ be a prime number and $m\in\mathbb N$ with $m\ge 2$, denote by $K$ the subgroup generated by $(p,0)$ and $(0,p)$ in $\G{p^m}$. If $A\subset\G{p^m}$ is a tiling complement of $K$, then 
$$A^{\Lsh}\cap \Z{A}=\emptyset.$$
\end{Lm}

\begin{proof}
Without loss of generality (by shifting, which does change $\Z{A}$) we may assume that $(0,0)\in A$. 

Let $a\in A$ be an arbitrary non-zero element, first we observe that 
\begin{equation} \label{EqSelf1}
|A\cap\bra{a^{\Lsh}}^{\perp}|\ge 2,
\end{equation}
since $A\cap\bra{a^{\Lsh}}^{\perp}$ contains at least $(0,0)$ and $(a^{\Lsh})^{\Rsh}=a$. 

Assume the contrary that $a^{\Lsh}\in\Z{A}$, then by \eqref{EqCount} of Theorem \ref{ThmCount} we have
\begin{equation} \label{EqSelf2}
|A\cap\bra{a^{\Lsh}}^{\perp}|=\frac{1}{p}\cdot |A\cap\bra{pa^{\Lsh}}^{\perp}|.
\end{equation}
Now there are a few cases concerning the subgroup  $\bra{pa^{\Lsh}}^{\perp}$, keeping in mind that $A$ is a tiling complement of $K$, we can assert that

(i) If $\bra{pa^{\Lsh}}^{\perp}$ is cyclic but not maximal, then it is contained in $K$, and $|A\cap \bra{pa^{\Lsh}}^{\perp}|=1$ (they intersect only at the identity); Alternatively if $\bra{pa^{\Lsh}}^{\perp}$ is cyclic and maximal, then $1<|A\cap \bra{pa^{\Lsh}}^{\perp}|\le p$;

(ii) If $\bra{pa^{\Lsh}}^{\perp}$ is not cyclic, then either it is contained in $K$, in which case we have $|A\cap \bra{pa^{\Lsh}}^{\perp}|=1$, or $K$ is contained in it with index $p$, in which case we have $|A\cap \bra{pa^{\Lsh}}^{\perp}|\le p$.

Therefore $|A\cap\bra{pa^{\Lsh}}^{\perp}|\le p$ must hold in all cases. Plugging it back into \eqref{EqSelf2} we get
$$|A\cap\bra{a^{\Lsh}}^{\perp}|\le 1,$$
which contradicts \eqref{EqSelf1}. 
\end{proof}

\begin{Cor} \label{CorDiff}
Let $p$ be a prime number and $m\in\mathbb N$ with $m\ge 2$, denote by $K$ the subgroup generated by $(p,0)$ and $(0,p)$ in $\G{p^m}$. If $A\subset\G{p^m}$ is a tiling complement of $K$, then 
$$\Delta(A^{\Lsh})\cap \Z{A}=\emptyset.$$
\end{Cor}

\begin{proof}
Assume the contrary that there exist distinct $a,b\in A$ such that $a^{\Lsh}-b^{\Lsh}\in \Delta(A^{\Lsh})\cap Z{A}$. Consider now the set
$$B=A-b,$$
then
$$B^{\Lsh}=A^{\Lsh}-b^{\Lsh}$$
contains $a^{\Lsh}-b^{\Lsh}$, and thus
\begin{equation} \label{EqDTemp}
a^{\Lsh}-b^{\Lsh}\in B^{\Lsh}\cap\Z{A}.
\end{equation}
Clearly $|B|=|A|$ and $\Delta B=\Delta A$ hold, hence $B$ is still a tiling complement of $K$. Consequently by Lemma \ref{LmSelf} we shall have 
$$B^{\Lsh} \cap\Z{B}=\emptyset.$$
On the other hand, it is also easy to see that $\Z{B}=\Z{A}$, therefore we shall have
$$B^{\Lsh} \cap\Z{B}=B^{\Lsh} \cap\Z{A}=\emptyset,$$
which contradicts \eqref{EqDTemp}. 
\end{proof}

If $m=2$, then Corollary \ref{CorDiff} indicates $\Delta A\cap\Z{A}=\emptyset$, which is the phenomenon described in the introduction. The condition that $A$ is a tiling complement of $K$ is sufficient but not necessary, in fact if $H$ is a subgroup, then from the relation (a counterpart of the Poisson summation formula)
$$\F{H}=|H|\cdot \mathbf 1_{H^{\perp}},$$
we see that the disjointness $\Delta(H^{\Lsh})\cap \Z{H}=\emptyset$ holds for all subgroup $H$. It also holds for tiling complements of proper subgroups in $\G{p}$.

\section*{Acknowldgement}
The author would like to thank Shilei Fan (CCNU) for his insightful inputs.

\end{document}